\documentclass[journal]{journal}

\usepackage{times}
\usepackage{amsmath,amssymb}
\usepackage{graphicx, wrapfig, subcaption, setspace, booktabs}%Include graphics package
\usepackage{setspace,fullpage}
\usepackage{upquote}% getting the right grave ` (and not ‘)!
\usepackage{float}%This is for the H in the include graphics command so pictures won't move
\usepackage{subcaption}%To add subfigures
\usepackage{braket}

\usepackage{mathdots} %for \iddots

\usepackage{geometry}
\usepackage{fancyhdr} % For custom headers
\usepackage{lastpage} % To determine the last page for the footer
\usepackage{extramarks} % For headers and footers
\usepackage[most]{tcolorbox} % For problem answer sections
\usepackage{xcolor} % For link coloring
\usepackage[hidelinks]{hyperref} % For URL links (no box or color name)
\usepackage{pdfpages}%to include pdf files
\usepackage{amsmath}

\usepackage{multicol}%to use columns

\usepackage{subfiles} % Best loaded last in the preamble

\usepackage{etoolbox}

\ifCLASSINFOpdf

\else

\fi

\begin{document}

\title{A Review on Higher Order Spline Techniques for Solving  Burgers Equation using B-Spline methods and Variation of B-Spline Techniques}

\author{Maryam~Khazaei Pool,~\IEEEmembership{}
        Lori~Lewis,~\IEEEmembership{}
       % <-this % stops a space
\thanks{M. Khazaei Pool is with the Department
of Electrical and Computer Engineering, University of California Merced,
CA, USA e-mail: (mkhazaeipool@ucmerced.edu}% <-this % stops a space
\thanks{Lori~Lewis is graduated master student with University of California Merced, CA, USA.}% <-this % stops a space
\thanks{Manuscript received June 30, 2022}}

\maketitle
\thispagestyle{empty}

\begin{abstract}
%\boldmath

This is a summary of articles based on higher order B-splines methods and the variation of B-spline methods such as Quadratic B-spline Finite Elements Method, Exponential Cubic B-Spline Method Septic B-spline Technique, Quintic B-spline Galerkin Method, and B-spline Galerkin Method based on the Quadratic B-spline Galerkin method (QBGM) and Cubic B-spline Galerkin method (CBGM). In this paper we study the B-spline methods and variations of B-spline techniques to find a numerical solution to the Burgers' equation. A set of fundamental definitions including  Burgers equation, spline functions, and B-spline functions are provided. For each method, the main technique is discussed as well as the discretization and stability analysis. A summary of the numerical results is provided and the efficiency of each method presented is discussed. A general conclusion is provided where we look at a comparison between the computational results of all the presented schemes. We describe the effectiveness and advantages of these methods.
\end{abstract}


\IEEEpeerreviewmaketitle

\section{Introduction}
Due to the great usefulness of spline functions in applications, scientists have used spline functions for various applications. Spline functions have applications in various fields such as applied mathematics and engineering. Spline methods are often used when solving Ordinary Differential Equations (ODE's) and Partial Differential Equations (PDE's).  B-spline \cite{e} methods have played an important role in computational mathematics, mathematical physics and mechanics. Geyikli and Gazi Karakoc applied septic B-spline collocation method for the numerical solution of the modified equal width wave equation \cite{f}. Parcha Kalyani, Mihretu Nigatu Lemma \cite{g} use a ninth degree spline function as well as an eighth degree spline to solve a seventh order boundary value problem. Authors obtained an approximate solution that very closely matches the exact solutions. Rashidinia and Khazaei solved a fifth degree and eight order boundary value problems with eight degree B-Spline \cite{i}. A low absolute error is obtained for their results which indicates that the presented numerical method is effective for solving high order boundary value problems. From these authors we get an understanding that spline methods produce solutions that are highly accurate.\\

\noindent Many efforts have been enhanced to evaluate the numerical solution of Burgers' equation in the past few years. The analytic solution of a two dimensional coupled Burgers' equation was first given by Fletcher \cite{k} using the Hopf-Cole transformation. A variety of studies have been developed for the various forms of nonlinear PDE's, as model problems in fluid dynamical systems \cite{m}, \cite{n}. Spline methods are commonly used to approximate a solution to Burgers' equation. A survey of higher order splines for boundary value problems by Srivastava can be found which gives a summary \cite{survey} of higher order spline methods.\\

\noindent In order to study the interactions of multi-shocks in thin viscoelastic tube filled,  S.Akter, M.G.Hafez \cite{npp} presented an analytic wave solutions of beta space fractional Burgers equation. This presented work investigates the single and overtaking collision of multi-shock wave excitations having space fractional evolution in a thin viscoelastic tube filled with incompressible inviscid fluid. The computational wave and numerical solutions of the Atangana conformable derivative $(1+3)$-Zakharov-Kuznetsov (ZK) equation with power-law nonlinearity are investigated via the modified Khater method and septic-B-spline scheme by M. A. Khater, Y. M Chu \cite{ppp}. This model is formulated and derived by employing the well-known reductive perturbation method.\\

\noindent The Burgers equation involves both non-linear propagation effects and diffusive effects. This equation is similar to the Navier–Stokes equation without the pressure term. Therefore, it is a simpler model to analyze fluid turbulence \cite{xcv}. For this paper we consider variations of B-spline methods including Quadratic B-Spline Finite Elements Method, Exponential-Cubic B-Spline Method and B-spline Galerkin methods for numerical solutions of Burgers equation. A summary of different methods for solving burgers' equation is presented. Using B-spline functions in different methods demonstrates efficient techniques in terms of time and cost in each approach. High accuracy solutions and stability and convergence analysis of the methods show efficiency and effectively of the discussed techniques. In Section \ref{sec:2.1}, we provide fundamental definitions of  Burgers equation, spline, and B-spline functions. Section \ref{sec:3} presents the quadratic B-spline finite element method which includes some remarks about the method and its advantages. In section \ref{sec:4}, the exponential-cubic B-spline method is summarized. Section \ref{sec:5} and \ref{sec:6} present the quintic B-spline Galerkin method and the septic B-spline techniques respectively. In section \ref{sec:7} we look at the B-spline Galerkin methods to find a numeric solution to Burgers' equation by considering two time splitting techniques. The final section, \ref{conclusion}, includes the conclusion and final thoughts about these methods.

\section{Definitions}
\label{sec:2.1}
\subsection{Burgers' Equation}
\noindent Burgers' equation is a fundamental partial differential equation \cite{BurgersEqn}. Here we will use the following form 

\begin{equation}
	u_t+uu_x=\lambda u_{xx}
	\label{BurgersEqn}
\end{equation}

% To reference Burgers' Eqn copy and paste the following (\ref{BurgersEqn})

\noindent where $\lambda > 0$ is a constant which is known as the diffusion coefficient,  and $u$ is an arbitrary function. \\

\noindent For each paper we describe the initial conditions
and boundary conditions that are considered for each method. 

%-------------------------------------
\subsection{Spline Functions}
\label{sec:2.2}
A spline is a piece-wise polynomial function defined in region $D=[a,b]$, such that there exists a decomposition of $D$ into sub-regions. In each sub-region of $D$, the function is a polynomial of some degree $k$.
The term "spline" is used to refer to a wide class of functions that are used in applications requiring data interpolation or smoothing. A function $S(x)$ is a spline of degree $k$ on $[a, b]$ if
\newline
$$ S\in C^{k-1} [a,b] $$
\\
$$ a=t_{0}<t_{1}<...< t_{n}=b $$
and

\begin{equation}
S(x)= 
\begin{cases}
S_{0}(x),&  t_{0} < x < t_{1}\\
S_{1}(x),              & t_{1} < x < t_{2}
\\
\hspace{5mm}\vdots &   \hspace{5mm}\vdots
\\
S_{n-1}(x),              & t_{n-1} < x < t_{n}
\end{cases}
\end{equation}

\noindent where $S_{i}(x)\in P^{k}, \quad i=0, 1, ..., n-1$.

%-------------------------------------
\subsection{B-Spline}
\label{sec:2}
The B-spline is defined as a basis function of degree $k$ which is denoted by $\varphi_{i}^{k}(x),$ where $i \in \mathbb{Z},$. In order to define B-spline basis functions, we need to define the degree of these basis functions, $p$. The i-th B-spline basis function of degree $p$ is written as $N_{i,p(u)}$ and it is defined as follows:

\begin{equation}
N_{i,p(u)}= 
\begin{cases}
1,& \text{if } u_{i+1} > u\geq u_{i}\\
0,              & \text{otherwise}
\end{cases}
\end{equation}

\begin{align}
N_{i,p(u)}=\frac{u-u_{i}}{u_{i+p}-u_{i}} N_{i,p-1}(u) + \\
\frac{u_{i+p+1} -u}{u_{i+p+1}-u_{i+1}} N_{i+1,p-1}(u) \nonumber
\end{align}

% where $w_{1}$ and $w_{2}$ 
\noindent The equation above is known as the Cox-de Boor recursion formula%CITATION for this
. To understand how the formula works, we can start by assuming the degree is zero (i.e., $p = 0$) then all the basis functions are considered as step functions. That is, basis function $N_{i,0}(u)$ is 1, if u is in the i-th knot span $[ui, u_{i+1})$. Here, we denote the B-spline of degree $k$ by $B_{i}^{k}(x),$ where $i$ is an  element in $\mathbb{Z}$ with the following properties \cite{asa}:\\

\noindent 1. $N_{i,p}(u)$ is a degree $p$ polynomial in $u$.
\\
2. Non-negativity: for all $i, p$ and $u$, $N_{i,p}(u)$ is non-negative.
\\
3. Local support: $N_{i,p}(u)$ is a nonzero polynomial on  $[u_i, u_{i+p+1}].$
\\
4. At most $p+1$ degree of the basis functions are nonzero
on any span $[ui, u_{i+p+1}]$,\\
\indent namely: $N_{i-p,p}(u), N_{i-p+1,p}(u), N_{i-p+2,p}(u),...,\\ N_{i,p}(u).$ This property shows that\\
\indent the following basis functions are nonzero on  $[u_i, u_{i+p+1}],$
$$
N_{i-p,p}(u), N_{i-p+1,p}(u), N_{i-p+2,p}(u), ..., N_{i,p}(u).
$$
\\
5. Partition of unity: The sum of all nonzero degree $p$, basis functions on span\\
\indent $[u_i, u_{i+p+1}]$ are 1 which states that the sum of these $p+1$ basis functions is $1$.
\\
\\
\\
An alternative approach to drive the B-Spline relations:
\\
Here we consider equally-spaced knots of a partition $\pi : a = x_0 < x_1 < . . . < x_n$ on $[a, b].$ This will be an alternative approach for deriving the B-splines which are more applicable with respect to the recurrence relation for the formulations of B-splines of higher degrees. Firstly, we recall that the $k-th$ forward difference $f(x_0)$ of a given function $f(x)$ at $x_0$ is defined recursively by \cite{c} and \cite{d} and is given as follows:

\begin{equation}
\nabla f(x_0) =\nabla(x_{1}) - \nabla f(x_{0}), 
\end{equation}

\noindent Definition: The function $(x-t)^{m}_{+}$, details given in
%Explain more about this citation 
\cite{i},
%where do you get this definition from? What is the definition off?
\begin{equation}
(x-t)^{m}_{+}=
\begin{cases}
(x-t)^{m} & 	x\leq t_{0}\\
0  &   x<t
\end{cases}
\end{equation}

\noindent It is clear that $(x-t)^{m}_{+}$ is $(m-1)$ times continuously differentiable with respect to $t$ and $x$.
The B-spline of order $m$ is defined as follows:

\begin{align}
B^{m}_{i}(t) = \dfrac{1}{h^{m}}\sum^{m+1}_{j=0} \begin{pmatrix} 
m+1   \\ 
j   \\     
\end{pmatrix} (-1)^{m+1-j}(x_{i-2+j}-t)^{m}_{+} \nonumber \\  =\dfrac{1}{h^{m}}\nabla^{m+1}(x_{i-2}-t)^{m}_{+} 
\end{align}

\noindent Hence, we can obtain the B-spline of various orders by taking various values of $m$.\\

\noindent Let $m = 1$ so that

\begin{align}
\dfrac{1}{h^{1}}\nabla^{2}(x_{i-2}-t)^{1}_{+} =\dfrac{1}{h^{1}}[(x_{i-2}-t)_{+} \nonumber \\ -2 (x_{i-1}-t)_{+} +(x_{i}-t)_{+}]
\end{align}
and 
\begin{equation} 
B^{m}_{i}(t) =
\begin{cases}
(x_{i-2}-t) -2 (x_{i-1}-t) & 	x_{i-2} < t \leq x_{i-1}\\
(x_{i}-t)  &   x_{i-1} < t \leq x_{i} \\
0  &   otherwise \\
\end{cases}
\end{equation}

\noindent By considering different values for $m$, different degree of B-Spline can be obtained, including septic B-Spline. 

%-------------------------------------
\section{Quadratic B-spline Finite Element Method}
\label{sec:3}
		For the quadratic B-spline finite element method presented in \cite{QuadBsplineMain}, the collocation method and a central difference with respect to time is used. This method is used to find a solution to the Burgers' equation (1) 
	with homogeneous boundary conditions that are
	\begin{equation}
		u(a,t)=u(b,t)=0,
	\end{equation}

	\noindent Since the finite element method  is used the region is partitioned into $N$ finite elements with equal length $h$ and knots $x_i$ are used such that $a=x_0<x_1<...<x_{N-1}=b$. The quadratic B-spline properties are defined as follows for $B_m$,
	
	\begin{align}
		B_m(x)=\frac{1}{h^2}\left\{
		\begin{array}{ll}
			(x-x_{m-1})^2, \\
			\quad  \text{if } x_{m-1}\leq x \leq x_m \\
			2h^2-(x_{m+1}-x)^2-(x-x_m)^2, \\
			\quad \text{if } x_m\leq x\leq x_{m+1} \\
			(x_{m+2}-x_m)^2, \\
			\quad \text{if } x_{m+1}\leq x \leq x_{m+2}\\
			0, \text{otherwise}
		\end{array}
		\right.
	\end{align}
	
	\noindent The goal is to approximate $u(x,t)$ of the form,
	\begin{equation}
		u(x,t)=\Sigma_m \xi_m(t)B_m(x)
	\end{equation}
	
	\noindent Here $\xi_m$ is given as a time dependent quantity and the numerical solution $u(x,t)$ is given in mid knots such as $y_m=(x_m+x_{m+1})/2$. The values of $u$ and the principal derivatives are calculated from the quadratic B-spline definitions,
	
	\begin{align}
		u&=u(x_m)=\frac{1}{4}\xi_{m-1}+\frac{3}{2}\xi_m+\frac{1}{4}\xi_{m-1} \nonumber \\
		u'_m&=u'(x_m)=\frac{-1}{h}\xi_{m-1}+\frac{1}{h}\xi_{m+1} \nonumber \\
		u''_m&=u''(x_m)=\frac{2}{h^2}\xi_{m-1}-\frac{4}{h^2}\xi_m+\frac{2}{h^2}\xi_{m+1}
	\end{align}
	
	\noindent To implement the collocation points,  mid knots are used to evaluate $u$ so that the following equation is obtained 
	\begin{align}
	    &\frac{1}{4}\xi_{m-1}^*+\frac{3}{2}\xi_m^*+\frac{1}{4}\xi_{m-1}^*+ \nonumber\\
		&\left(\frac{1}{4}\xi_{m-1}+\frac{3}{2}\xi_m+\frac{1}{4}\xi_{m-1}\right) \nonumber\\
		&\left(\frac{-1}{h}\xi_{m-1}+\frac{1}{h}\xi_{m+1}\right) \nonumber\\
	    &=\nu\left(\frac{2}{h^2}\xi_{m-1}-\frac{4}{h^2}\xi_m+\frac{2}{h^2} xi_{m+1}\right) \nonumber
	\end{align}
	
	\noindent Here the $*$ represents the differentiation with respect to time. Interpolating between $n$ and $n+1$ and using the central difference operator for time, $\xi$ , gives
	\begin{equation}
		\xi=\frac{\xi^{n+1}+\xi^n}{2}, \quad \xi^*=\frac{\xi^{n+1}-\xi^{n-1}}{2\Delta t}
	\end{equation}
	
	\noindent where $\xi^n$ are parameters at the time $n\Delta t$. A system is then obtained and written as 
	\begin{equation}
		A(\xi)\xi^{n+1}=B(\xi)\xi^n+C\xi^{n-1}
	\end{equation}
	
	\noindent here $A(\xi)$, $B(\xi)$ are tridiagonal matrices with two initial time levels. The exact solution at $t=t_0$ and $t=t_0+\Delta t$ is used to obtain the initial conditions.\\
	
	\noindent Stability Analysis: Since the finite element method is used, which is explained in \cite{QuadBspline2}, an investigation of the stability of the numerical scheme is required. The Von-Neummann method is used to find the stability which is defined as
	\begin{equation}
		\xi_m^n=\hat{\epsilon}^n e^{imkh}
	\end{equation}
	Here $k$ is the mode number and $h$ is the element size.
 The following equations are then solved and the roots $g_1$ and $g_2$ are found as
	\begin{equation}
		g_1=-1,\quad g_2=\frac{3+\cos \theta}{3+\cos \theta+4\alpha\cos \theta+2\beta}
	\end{equation}
	
	where $\alpha=\Delta z/4-2\Delta t\nu/h^2$, $\beta=3\Delta t z/2+4\Delta tv/h^2.$
\noindent	Modules of growth  are taken to obtain $|g|\leq1$, which means that the scheme is unconditionally stable.\\
	
\noindent Remarks:	The $L_2$ and $L_\infty$ error norms are used to compare the analytical and numerical solutions. Different examples with different initial conditions are analyzed. Here we focus on the initial condition of an exponential form for which the exact solution is known to be 
\begin{equation}
	u(x,t)=\frac{x/t}{1+\sqrt{(t/t_0)e^{x^2/(4\nu)}}}, \quad t\geq 1.
	\label{ExponentialExample}
\end{equation}

% to reference the top equation use this -> (\ref{ExponentialExample})

\noindent The results are compared with the exact solution for different $h$, and $k$ values. The algorithm is compared with the exact solutions and it is shown that the method produces accurate results for small viscosity values. The quadratic B-spline method  is easy to implement and can be generalized with higher order spline methods. Because of the flexibility and accuracy of this method the quadratic B-spline finite element method is advantageous when finding a solution to the Burgers' equation.

%---------------- Next paper ---------

\section{Exponential-Cubic B-Spline Method}
\label{sec:4}
The exponential cubic B-spline functions  are used to set up the collocation method to solve Burgers' Equation by \cite{ExponentialCubic}. The initial conditions and boundary conditions considered are as follows;\\

\noindent the initial condition 
\begin{equation}
	u(x,0)=f(x), \quad a\leq x\leq b
\end{equation}
and the boundary conditions 
\begin{equation}
	u(a,t)=\alpha_1, \quad u(b,t)=\alpha_2
\end{equation}
Here $\alpha_1,\quad \alpha_2$ are constants, $u = u(x,t)$ is
a sufficiently differentiable unknown function and $f(x)$ which is bounded. \\

\noindent Exponential Cubic B-spline Collocation Method:
For this method the nodes are equally distributed for the domain so that 
\begin{equation}
	\pi: \quad a=x_0<x_1< \hdots <x_N=b
\end{equation}
and a mesh with spacing $h=(b-a)/N$ is used. The exponential cubic b-splines, $B_i(x)$, at the points of $\pi$ are defined as

\begin{equation}
	B_i(x) = \left\{
	\begin{array}{ll}
		b_2((x_{i-2}-x)-\frac{1}{p}(\sinh(p(x_{i-2}-x)))) & [x_{i-2},x_{i-1}], \\
		a_1+b_1(x_i-x)+c_1 \exp(p(x_i-x))+d_1\exp(-p(x_i-x))& [x_{i-1},x_i],\\
		a_1+b_1(x-x_i)+c_1 \exp(p(x-x_i))+d_1\exp(-p(x-x_i))& [x_i,x_{i+1}],\\
		b_2((x-x_{i+2})-\frac{1}{p}(\sinh(p(x-x_{i+2})))) & [x_{i+1},x_{i+2}], \\
		0 & \text{otherwise}
	\end{array}
	\right.
\end{equation}

where
\begin{align}
	a_1=& \frac{phc}{pch -s}, \quad b_1=\frac{p}{2}[\frac{c(c-1)+s^2}{(pch-s)(1-c)}], \nonumber
\end{align}
and

\begin{align}
	c_1=& \frac{1}{4}[\frac{\exp(-ph)(1-c)+s(\exp(-ph)-1)}{(pch-s)(1-c)}]
\end{align}

\begin{align}
	d_1=&\frac{1}{4}[\frac{\exp(ph)(c-1)+s(\exp(ph)-1)}{(phc-s)(1-c)}]
\end{align}

and $c=\cosh(ph)$, $s=\sinh(ph)$, $p$ is a free parameter. A basis is formed for the functions over the interval. Each basis function $B_i(x)$ is twice continuously differentiable and $B_i$, $B_i'(x)$ and $B_i''(x)$ \cite{ExponentialCubic}. To approximate the unknown $u$, $u_N$ is used which is in the form of
\begin{equation}
	u_N(x,t)=\sum_{i=-1}^{N+1}\delta_iB_i(x)
\end{equation}
Here $\delta_i$ is a time dependent parameter. The first and second derivatives are calculated at Knots, $x_i$, and the Crank-Nicolson scheme is used to discretize time variables of the initial conditions are modified known $u$ in the Burgers’ equation. After some substitution the initial conditions at the boundaries are obtained and used to find an approximation to the Burgers' equation.\\

\noindent Remarks: In a similar approach with the previous sections, the discrete $L_2$ and $L_\infty$ error norms are used to compare the analytical and numerical solutions. Similar to the quadratic B - spline finite element method an example considered is a particular solution to Burgers equation which has the following initial condition
\begin{equation}
	u(x,1)=exp\left(\frac{1}{8\lambda}\right), \quad 0\leq x \leq 1,
	\label{initalCondition}
\end{equation}
with boundary conditions $u(0,t)=0$ and $u(1,t)=1$. The reason why this example is chosen is because the analytical solution is known to be (\ref{ExponentialExample}).
It is mentioned that the solution to this specific example will be successful for a small $\lambda$ which results in a steep shock solution. The propagation of the shock is studied with $\lambda=0.005$ and $\lambda=0.0005$. The results are compared with other papers \cite{CubicBsplines} and \cite{GalerkinBspline}. The exponential cubic B-spline collocation method provides better results than the cubic B-spline collocation method and the B-spline Galerkin finite element method. Note that the cost of the cubic B-spline Galerkin method is higher than the exponential cubic B-spline method. Over all the test runs of the exponential cubic B-spline method had the best results for finding the free parameter $p=1$.

%---------------- Next paper ---------
\section{Quintic B-Spline Galerkin Method}
\label{sec:5}
Now we look into the quintic B-splines method to find solutions of a time-split Burgers equation over finite intervals with the help of \cite{QuinticBspline}. The Burgers equation studied has the following boundary conditions 
\begin{align*}
	u(a,t)&=\alpha_1, \quad  \quad u(b,t)=\alpha_2,\\
	u_x(a,t)&=0,  \quad  \quad u_x(b,t)=\alpha_2, \quad \quad t\in (0,T],\\
	u_{xx}(a,t)&=0,  \quad  \quad u_{xx}(b,t)=\alpha_2
\end{align*}

\noindent The paper explores both solutions for the Burgers and time-split Burgers equation, but here we focus on the solution of the Burgers equation. \\

%-------------------------------------

\noindent Quintic B-spline Galerkin method:
The method begins by applying the Galerkin technique\cite{GalerkinBspline}  to add the weighted functions and a mesh $a=x_0<x_1< \hdots <x_N=b$ as an uniform partition. Here the knots $x_m$ and $h=x_m-x_{m-1}$, $=1,...,N$.\\
$Q_m(x), \quad m=-2,...,N+2$ is given as

\begin{equation*}
    \text{For the following formula, let } w_i = (x-x_{m-i})^5
\end{equation*}

\begin{equation*}
	Q_m(x) = \frac{1}{h^5}\left\{
	\begin{array}{l}
		w_3, \nonumber \\ \quad \text{if } [x_{m-3},x_{m-2}], \\
		w_3-6w_2, \nonumber \\ \quad \text{if } [x_{m-2},x_{m-1}],\\
		w_3-6w_2+15w_1, \nonumber \\ \quad \text{if } [x_{m-1},x_m],\\
		w_3-6w_2+15w_1-20w_0, \nonumber \\ \quad \text{if } [x_m,x_{m+1}], \\
		w_3-6w_2+15w_1-20w_0+15w_{-1}, \nonumber \\ \quad \text{if } [x_{m+1},x_{m+2}], \\
		w_3-6w_2+15w_1-20w_0+15w_{-1}-6w_{-2}, \nonumber \\ \quad \text{if } [x_{m+2},x_{m+3}], \\
		0, \text{otherwise}
	\end{array}
	\right.
\end{equation*}

\noindent The quintic B-splines with knots $x_m$, $m=-5,...,N+5$ and bases on $a\leq x\leq b$ with the global approximation defined as
\begin{equation}
	u_N(x,t)=\sum_{m=-2}^{N+2} \delta _m(t)Q_m(x)
\end{equation}	
where $\delta _m$ is a time dependent parameter. 
A local coordinate system is used. It is defined by the mapping the relationship $\xi=x-x_m$ to obtain a global and local coordinates using transformations of finite element $[x_m,x_{m+1}]$ into interval $[0,h]$.

The approximation is reduced over the element $[x_m,x_{m+1}]$ as follows
\begin{equation}
	u_N^e=u(\xi,t)=\sum_{i=m-2}^{m+3}\delta_i(t)Q_i(\xi),
\end{equation}

\noindent where $\delta_i$, $i=m-2,...,m+3$ are element parameters. The weighted function with the quantic B-spline  and  $U_N^e$  over the element $[0,h]$ gives a weighted function as a matrix form
\begin{equation}
	A^e\dot{\delta^e}+(\delta ^e)^TL^e\delta^e-\lambda D^e\delta^e
\end{equation}
It is noted that the matrices $A$, $D$ are 6 by 6 and the matrix $L$ is a 6 by 6 by 6 which $\delta^e$ is defined as follows:\\
\begin{equation}
\delta^e=(\delta_{m-2},\delta_{m-1},\delta_m, \delta_{m+1},\delta_{m+2}, \delta_{m+3})
\end{equation}

\noindent The matrix $L$ is organized as the following
\begin{equation}
	B_{i,j}=\sum_{k-m-2}^{m+3}L_{ijk}\delta_k
\end{equation}

\noindent After combining all element matrices as a system of nonlinear ordinary differential equations with a global parameter $\delta$, the following equation is obtained:
\begin{equation}
A\delta^{0}+(B-\nu D) \delta=0
\end{equation}

\noindent using the Crank-Nicolson discretization formula for the vector of element
parameter $\delta$ and the finite difference equation for the time derivatives parameters
$\delta^{0}$,  a nonlinear recurrence relation for the time parameters $\delta$ with the following form is obtained:

\begin{align}
(2A+\nabla(t)B -\nu \nabla(t)D) \delta^{n+1}=(2A-\nabla(t)B \nonumber\\
+\nu \nabla(t)D) \delta^{n}
\end{align}

\noindent Boundary conditions at the left end of the region and and at the right end of the region are applied as well as some terms of the global parameter so that a solvable system made of $N+5$ equations and $N+5$ unknown parameters is obtained.  An 11-banded matrix system at every time step is solved.
To start the iteration of the recurrence relation of the above system, the initial parameter vector $\delta^{0}$ is obtained and time evaluation of $u_{N}$ can be evaluated from the time evolution of the vector $\delta^{n}$, which is found by solving the recurrence relation of the system above.\\

%-------------------------------------
\noindent Remarks:
The $L_2$ and $L_\infty$ error norms are used to compare the analytical and numerical solutions. Just like other methods here we focus on the first test problem with an exponential initial condition (\ref{initalCondition}) which the analytical solution to Burgers' equation is known (\ref{ExponentialExample}).
Here the parameters used are $\lambda = 0.005$, $h = 0.005$ and $\Delta t = 0.01$ over the problem domain $[0, 1]$. This method shows accuracy in the $L_2$ norm. When compared to the cubic spline methods \cite{CubicBsplines}, \cite{ExponentialCubic} it is noted that the method shows small improvement for the $L_2$ and $L_\infty$ norms.\\

\noindent The quintic B-spline  method provides provide high accuracy results for finding the solution of Burgers' equation. The time splitting does not affect the method. It was concluded that the method is efficient and reliable.

%-------------------- Next Paper ----------------
\section{Septic B-spline Techniques to Solve Burgers Equation}
\label{sec:6}
 The septic B-spline method over finite elements \cite{p}  is used to obtain a numerical solutions to the nonlinear Burgers' equation by considering  \cite{a}. The papers focus on obtaining a solution for Burgers equation (\ref{BurgersEqn}) with the initial condition,
\begin{equation}
u(x,0)=f(x)
\end{equation}

and the following boundary conditions
\begin{equation}
u(a,t)=\alpha_{1},\quad u(b,t)=\alpha_{2},
\end{equation}

and 

\begin{equation}
\quad u_{x}(a,t)=u_{x}(b,t)=0,
\end{equation}

\begin{equation}
u_{xx}(a,t)=u_{xx}(b,t)=0,
\end{equation}

\begin{equation}
u_{xxx}(a,t)=u_{xxx}(b,t)=0,
\end{equation}

\noindent where $\alpha_{1}, \alpha_{2}$ are constants as the problem need, $u = u(x, t)$ is a sufficiently-often differentiable function, and $f(x)$ is bounded.\\

\noindent The authors \cite{a} discussed about finding approximate solution $u_{N}(x, t)$ which satisfy the following conditions:
(a) It must agree with the initial condition $u(x,0)$ at the knots $x_{j}$.
(b) The first, second and the third derivatives of the approximate initial condition agree with the exact initial conditions at both ends of the range $[a, b]$.
The approximate solution to $u(x,t)$ is in the form of a collocation method:
\begin{equation}
u_{N}(x, t)=\sum_{i=-3}^{N+3} B_{i}(x_{j})\omega_{i}(t), \quad j=0, 1,..,N
\end{equation}

\noindent where $x_{i}(t)$ are time dependent quantities to be determined. Using the boundary conditions where $\varPhi_{i}(x_{j})$ is the values of the septic
B-spline function forms a basis for the functions defined over $[a,b]$, and all its first, second, and third derivatives vanish outside the interval $[x_{i-4},x_{i+4}]$. Here $\omega_{i}(t)$ is a time dependent quantity that is determined by using the given boundary conditions in the paper. Using values of the bases in the collocation method and its derivatives with knots at the shown points produces a matrix system that consists of $N + 1$ equations with $N + 7$ unknowns. This requires six additional constraints which are obtained from the boundary conditions. The presented matrix system equation has the following form:

\begin{equation} A(\omega^{n})\omega^{n+1}=B(\omega^{n})\omega^{n}+r  
\end{equation}

\noindent where the matrices $A(\omega^{n})$, and $B(\omega^{n})$ are septa-diagonal$(N + 1)\times (N + 1)$ matrices and $r$ is the $N + 1$ dimensional column vector.
The septa-diagonal algorithm is then used to solve the derived system to obtain a solution.\\

\noindent Stability Analysis: Stability analysis is done by using the Von-Neumann stability analysis for the linear system. An amplification factor $g$ is obtained for mode $k$ and results produce $|g|\leq 1$  which means that the linearized numerical scheme for the Burgers' equation is unconditionally stable.\\

\noindent Remarks: 
The numerical solutions of Burgers' equation and modified Burgers' equations are analyzed by computing the difference between the analytic and numerical solutions at each mesh point. The $L_{2}$ and $L_{\infty}$ norms are used for comparison of the results.
The results of the nonlinear Burgers' equation by Septic B-spline technique are given as follows:\\
\noindent (1): As the viscosity value $\lambda$ is increased the errors tend to increase, but for all the values of  used here, the errors are acceptable.

\noindent (2): It is discussed that as the time increases, the curve of numerical solution decays.

\noindent (3): The  numerical solutions obtained exhibit to maintain good accuracy compared with the exact solution, especially for small values of the viscosity parameter. Using a collocation method with the septic B-splines gives an accurate approximation, particularly for small values of the viscosity parameter.

%---------------- NEW PAPER --------------------------------
\section{B-spline Galerkin Methods for Numeric Solutions of the Burgers' Equation}
\label{sec:7}
Here we will look at one variation of B-Spline functions to solve Burgers' equation \cite{b}. A solution to Burgers equations is approximated using quadratic and cubic B-spline Galerkin finite element method. The Burgers' equation (\ref{BurgersEqn}) solved has the following initial condition,

\begin{equation}
u(x,0)=f(x)
\end{equation}

\noindent and boundary conditions
\begin{equation}
u(a,t)=\alpha_1,\quad u(b,t)=\alpha_2, \quad u_{x}(a,t)=u_{x}(b,t)=0
\end{equation}

\noindent where subscripts $x$ and $t$, where $t\in [0, T]$ denote differentiation.\\

\noindent A system of PDE's of the first order is obtained by splitting Burgers' equation in time as follows

\begin{equation}
u_{t}+2uu_{x}=0, \quad
u_{t}-2vu_{xx}=0
\end{equation}

\noindent Applying the Galerkin technique and considering weight functions, $w$, to the equations above lead to following integral equations,
\begin{equation}
\int ^{b}_{a} w( u_{t}+2uu_{x} dx)=0,  \quad
\int ^{b}_{a} w(u_{t}-2vu_{xx} dx)=0
\end{equation}

\noindent The first and second order smooth solutions are provided by using the quadratic and cubic B-splines functions as well as the Galerkin finite element method \cite{b}. \\

\noindent The Quadratic B-spline Galerkin method (QBGM):
Here a global approximation $u_N$ is written in terms of B-splines \cite{b} given in the following form :
\newline
\begin{center}
\begin{equation}
u_{N}(x, t)= \sum_{m=-1}^{N} \delta_{m}(t) Q_{m}(x),
\end{equation}
\end{center}

\noindent where $\delta_{m}$ is time dependent parameter which is specified from the quadratic Galerkin method. $Q_{m}(x)$ represents the quadratic B-splines at knots $x_{m}$. A basis is then formed over the interval $[a, b]$. First derivatives values vanish outside the interval. An interval $[x_{m}, x_{m+1}]$ includes three successive quadratic B-splines.
\newline
The finite elements are identified with intervals  $[x_{m},x_{m+1}]$ with nodes at
$x_{m}$ and $x_{m+1}$. This transforms the quadratic B-splines into element shape functions over the finite intervals $[0,h]$. A local coordinate system is used which is $\xi = x - x_{m}$, where $\xi \in [0, h]$.
\
\noindent Consider $\delta^{e} = (\delta_{m-1}, \delta_{m}, \delta_{m+1})$ which are known as element parameters and $Q^{e} = (Q_{m-1}, Q_{m}, Q_{m+1})$ are given as element shape functions. A system of algebraic equations is obtained by applying the Galerkin method, and considering both weight and approximate functions that are chosen as the quadratic B-spline shape functions:
\newline
\begin{equation}
(2A+ \nabla t B)\delta^{n+1/2}= (2A - \nabla t B)\delta^{n}
\end{equation}

\begin{equation}
(2A- v\nabla t D)\delta^{n+1}=(2A+ v\nabla t D)\delta^{n+1/2}
\end{equation}

\noindent The two pentadiagonal systems, shown above, consist of $(N + 2)$ equations of $(N + 2)$ unknown parameters
$
(\delta^{n}_{-1}, \delta^{n}_{0}, \delta^{n}_{N}). 
$
Applying the boundary conditions $u(a,x) = u(b,x) = 0$ at both ends of the interval and using Thomas algorithms, the solutions of the pentadiagonal matrix equations with the dimensions $N \times N$ are obtained.
After initial parameters $\delta^{0}_{m}$ are obtained with the help of the boundary and initial
conditions, time evolution of the parameters $\delta^{n}_{m}$ are computed using the recurrence relations between time steps. These are obtained alliteratively so that time evolution of the approximate solution $u_{N}$ could be determined \cite{b}.\\

\noindent Cubic B-spline Galerkin method (CBGM):
Cubic B-spline $Q_{m}, \quad m=-1, ..., N+1$ is defined at the knots $x_{m}$ and a basis is formed over $[a, b]$. An approximate solution to $u_{N}(x,t)$ is obtained using the cubic B-splines \cite{b} and element parameters $\delta_{m}$ which have the following form :
\newline

\begin{equation}
u_{N}(x, t)= \sum_{m=-1}^{N+1} \delta_{m}(t) Q_{m}(x).
\end{equation}

\noindent Using the above expression and the values of the cubic B-splines $Q_{m}$ at the knots $x_{m}$, the values of $u, u'$and $u"$ in terms of the element parameters are obtained and given as follows:
\newline
\begin{align}
 u_{m}&=\delta_{m-1}+4\delta_{m}+\delta_{m+1},\\ 
 u'_{m}&=\dfrac{3}{h}(\delta_{m+1}-\delta_{m-1}),\\
 u"_{m}&=\dfrac{6}{h^{2}}(\delta_{m-1}-2\delta_{m}+\delta_{m+1}),
\end{align}

\noindent where time dependent parameters, $\delta_{m}$, is determined from the cubic B-spline Galerkin method.\\

\noindent A mapping of a typical finite interval $[x_{m},x_{m+1}]$ to the interval $[0, h]$ is used with local coordinates $\xi$ to related to the global coordinates $ \xi=x-x_{m}, \xi\in [0, h]$. Using the given cubic B-spline shape functions $Q_{m-1}, Q_{m}, Q_{m+1}, Q_{m+2}$  in terms of the $\xi$ over the $[0, h] $ will cover a finite element $[x_{m}, x_{m+1}]$ which yields to a local approximation (trial solution) over the element and is given as follows,
\begin{align}
u^{e}_{N}=Q_{m-1}(\xi) \delta_{m-1}(t)+ \nonumber \\
Q_{m}(\xi) \delta_{m}(t)+Q_{m+1}(\xi) \delta_{m+1}(t)+ \nonumber \\
Q_{m+2}(\xi) \delta_{m+2}(t),\nonumber\\ 
\end{align}

where the element parameters are $\delta^{e}=(\delta_{m-1}, \delta_{m}, \delta_{m+1}  \delta_{m+2})$ and element shape functions are $Q^{e}= (Q_{m-1}, Q_{m}, Q_{m+1}, Q_{m+2})$. Substituting weight functions $W$ and $u$ by shape functions $Q_{m}$ and trial solution (67) into the main equations (58) yields to a matrix system of first order
ODE's which leads to a global matrix equation:
\begin{equation}
(2A+ \nabla t B)\delta^{n+1/2}= (2A - \nabla t B)\delta^{n}
\end{equation}

\noindent In a similar manner, by using interpolation of the parameters $\delta_{m}$ and it's time derivative between two time levels $n + 1/2$ and $n + 1$, an algebraic equation is obtained as follows:
\begin{equation}
(2A- \lambda\nabla t D)\delta^{n+1}=(2A+ \lambda\nabla t D)\delta^{n+1/2}
\end{equation}
The equations (68) and (69) consist of two recurrence relations for the time $(N + 3)$ equations
of the $(N + 3)$ unknown parameters. Applying boundary conditions produces a septa-diagonal systems which includes $(N + 1)$ unknown parameters in equations $(N + 1)$. Then the time evolution of the time parameter for the both schemes are obtained \cite{b}.\\

\noindent Remarks:
The numerical solution of Burgers' equation is discussed \cite{b} for three standard problems. The $L_{2}$ and $L_{1}$ error norms are used to measure the versatility and accuracy of the proposed methods as well as $|e_{1}|$ norm. The Galerkin method with both quadratic and cubic B-splines are presented as weight and trial functions which are used to obtained a solution to the time-split Burgers equation. \\

 \noindent The first example is about Shock-like solution of the Burgers' equation which is compared with the analytical solution (\ref{ExponentialExample}). The propagation of the shocks is shown to be  slightly smoother as time increases. A variety of boundary conditions are tested and the best result is obtained by selecting zero for initial conditions as $u(a, t)=0$ and $u(b, t)=0$. Both schemes show the same result for  the $L_{2}$ and $L_{1}$ error norms. From the results, the present calculation produced has a larger error compared to the schemes in which the split Burgers' equation approximation is not carried out.
\\

\noindent For the second example the Burgers equation is discussed with the following initial condition
\begin{equation}
u(x,0)= \sin(\varPi x), \quad x\in[0,1]
\end{equation}
and boundary condition
\begin{equation}
u(0,t)=u(1,t)=0, \quad 0\leq t.
\end{equation}
We observe the decay of sinusoidal disturbance. The parameters which are used are viscosity constant $m = 1$, time step $\Delta t = 0.00001$ and various space steps are considered. There is a good agreement between both numerical schemes and exact values. Numerical results for $\lambda = 10^{4}$ show a very sharp front near the left boundary at earlier times and as time increases. The sharpness and amplitude of the wave front then decays. These properties of the numerical solutions from the QBGM and CBGM produce a small error when comparing with the result obtained by Varog˘lu and Finn \cite{aa}, Kakuda and Tosaka \cite{bb}.  For an arbitrary initial data sets, the exact solutions of Burgers equation is presented as a quasi-linear parabolic PDE. Considering the fact that the analytical solutions of Burgers equations involve Fourier series solutions for a small viscosity constant $\nu$, the analytical solutions converge slowly.

%**************************************8
\section{Conclusion}
\label{conclusion}
\noindent We provide a summary of different methods for solving Burgers' equation which are shown to be efficient and effective. A summary for solving time-split Burgers' equation is presented by using quadratic and cubic B-spline Galerkin finite element techniques. Solving Quintic B-spline Galerkin Method results in an 11-banded sparse matrix system for every time step which is efficient time wise and cost wise. Two numerical algorithms based on Galerkin method with both quadratic and cubic B-splines as weight and trial functions are studied for the time-split Burgers' equation.  This technique produces a high accuracy solution for Burgers' equation. Moreover, by selecting suitable boundary conditions for the Galerkin method with both cubic and quadratic B-splines as an approximate function will produce a similar error. Having sparse and band matrices in a linear system for the septic B-spline function techniques \cite {a} is more efficient and cost less computationally. Stability analysis for the methods show that the methods are stable which are great to work with. Finally we present a comparison among the numerical results of all schemes and analytical values in all methods which maintain a good accuracy compared with the exact solutions. These methods are efficient and cost effective and are a great option for solving Burgers' equation.  
%\begin{acknowledgements}
%\end{acknowledgements}



\begin{thebibliography}{}

    \bibitem{CubicBsplines}  \.{I}. Da\u{g}, D. Irk and B. Saka: A numerical solution of the Burgers’ equation using cubic B-splines, \emph{Appl. Math. Comput.}, 163, 199-211 (2005)

	\bibitem{GalerkinBspline} \.{I}. Da\u{g}, B. Saka,  and A. Boz: B-spline Galerkin methods for numerical solutions of the Burgers’ equation, \emph{Appl. Math. Comput.}, 166, 506-522(2005)
	
	\bibitem{b}
	I. Dag, B. l. Saka, A. Boz: B-spline Galerkin methods for numerical solutions of the Burgers equation, \emph{Elsevier, Applied Mathematics and Computation}, 506–522(2005)

	\bibitem{ExponentialCubic} O. Ersoy and I. Dag and N. Adar: The Exponential Cubic B-spline Algorithm for Burgers's Equation,  \emph{ arXiv: Numerical Analysis}, (2016)
		
    \bibitem{k}
	C. Fletcher: Generating exact solutions of the two-dimensional Burgers equation, \emph{Int. Numer. Meth. Fluids}, 203, 213-216(1983)
		
	\bibitem{f}
	T. Geyikli, S B Gazi Karakoc: Septic B-Spline Collocation Method for the Numerical Solution of the
	Modified Equal Width Wave Equation, \emph{Inonu University, Malatya, Turkey.}(2011)
		
	\bibitem{bb}
	K. Kakuda, N. Tosaka: The generalized boundary element approach to Burgers equation, \emph{Int.
	J. Numer. Meth. Eng.}, 245–261(1990)
		
	\bibitem{e}
	P. K. Srivastava, Study of differential equations with their polynomial and non-polynomial spline based approximation, \emph{Acta Tehnica Corviniensis Bulletin of Engineering Tome VII}, 2067 3809(2014)

	\bibitem{QuadBspline2} S. Kutluay, A. Esen and I. Dag: Numerical solutions of the Burgers’ equation by the least squares quadratic B-spline finite element method,  vol. 167, pp. 21-33.   \emph{J. Comput. Appl. Math.,} 80, 931–938(2004)
		
    \bibitem{ee}	Bickely, W.G.: Piecewise Cubic Interpolation and Two-Point Boundary Value Problem, Computer Journal, pp. 202-208(1968)
		
	\bibitem{n}
	W. Liu: A asymptotic behavior of solutions of time-delayed Burgers equation, \emph{ Discrete and Continuous Dynamic Systems,Series B}, vol. 2, no. 1, pp. 47–56(2002)
		
	\bibitem{p}
	H. Nguyen, J. Reynen: A space time finite element approach to Burgers equation,  \emph{E.
	Hinton et al. (Eds.)}, Numerical Methods for Non-linear Problems, vol. 3, Pineridge Press, pp. 718–728(1987)
	
	\bibitem{g}
	K. Parcha, N. L. Mihretu: Solutions of Seventh Order Boundary Value Problems Using Ninth Degree Spline Functions and Comparison with Eighth Degree Spline Solutions, \emph{Journal of Applied Mathematics and Physics}, vol. 4, no.2, 249-261(2016)
			
	\bibitem{d}
	P. M Printer: Splines and variational Methods, \emph{Colorado State University, Wiley Classics Edition published}, 57, 421-421 (1975)
				
	\bibitem{a} 
	M. A. Ramadan, T. S. El-Danaf, F. E.I. Abd Alaal: A numerical solution of the Burgers equation using septic B-splines, \emph{Chaos, Solitons and Fractals}, 795–804(2005)
	
	\bibitem{i}
	J. Rashidinia, M. Khazaei, and H. Nikmarvani, Spline collocation method for solution of higher order linear boundary value problems, \emph{TWMS J. Pure Appl. Math.}, 6, 38-47(2015)
	
	\bibitem{az}
	Bateman H.: Some recent researches on the motion of fluids, Monthly Weather Review, 163–70(1915)
		
	\bibitem{bz}
	Burger JM.: A mathematical model illustrating the theory of turbulence, Advanced in Applied Mechanic I. New York, 171–99(1948)
		
	\bibitem{c}
	J. Rashidinia, and Sh. Sharifi, Survey of B-spline functions to approximate the solution of mathematical problems, \emph{Iran University of Science and Technology.}, (2011)
	
	\bibitem{QuadBsplineMain} K. R. Raslan, A collocation solution for Burgers equation using quadratic B-spline finite elements,  \emph{ International Journal of Computer Mathematics,}, 80, 931–938(2003)
	
	\bibitem{QuinticBspline} B. Saka and \.{I}. Da\u{g} and A. Boz, Quintic B-Spline Galerkin Method for Numerical Solutions of the Burgers' Equation(2004)
	
	\bibitem{m}
	T. E. SAYED, A. E. Danaf, Numerical solution of the Korteweg–Vries Burgers equation by quintic spline method, \emph{Studia Unv, Babes––Bolyai Math}, 41–55(2002)
	
	\bibitem{BurgersEqn}
	M. Shearer and R. Levy, Partial Differential Equations Introduction to Theory and Applications, \emph{Princeton university press}, 175-176(2015)
	
	\bibitem{j}
	H. S. Shukla , M. Tamsir, V. K. Srivastava, and J. Kumar:(2014), Numerical solution of two dimensional coupled viscous Burgers equation using modified cubic B-spline differential  quadrature method, \emph{AIP (American Institute of Physics (United States)) Advanced}, 4, (2015) 
	
	\bibitem{survey} 
	P. K. Srivastava: Application of higher order splines for boundary value problems. \emph{International Journal of Mathematical, Computational, Statistical, Natural and Physical Engineering 9.2}, 115-122(2015)
	
	\bibitem{aa}
	E. Varog lu, W.D.L. Finn, Space–time finite elements incorporating characteristics for the Burgers equation, \emph{Int. J. Numer. Meth. Eng}, pp. 171–184(1980)
	
	\bibitem{asa}
	M. Khazaei, Y. Karamipour, Numerical Solution of The Seventh Order Boundary Value Problems using B-spline Method, arXiv:2109.06030v1 [math.NA], arXiv:2109.06030 [math.NA]
	
	
	
    \bibitem{ppp}
    M. A. Khater, Y. M Chu, R. A. M. Attia, M. Inc, and D. Lu, On the Analytical and Numerical Solutions in the Quantum Magnetoplasmas: The Atangana Conformable Derivative (1+3)-ZK Equation with Power-Law Nonlinearity, https://doi.org/10.1155/2020/5809289, \emph{Hindawi, Advances in Mathematical Physics}, (2020) 
    
    \bibitem{npp}
    S.Akter, M.G.Hafez, Yu-Ming Chu and M.D.Hossain, Analytic wave solutions of beta space fractional Burgers equation to study the interactions of multi-shocks in thin viscoelastic tube filled, \emph{Science Direct}, pp. 877-887(2021)
	
	\bibitem{xcv}
	H. Ramos, A. Kaur and V. Kanwar, Using a cubic B-spline method in conjunction with a one-step optimized hybrid block approach to solve nonlinear partial differential equations, \emph{Computational and Applied Mathematics}, (2022) 
\end{thebibliography}


\end{document}